\documentclass[12pt]{article}
\usepackage{amsmath,amssymb}

\def\seq#1#2#3{#1_{#2},\,\ldots,#1_{#3}}
\def\w{\widetilde}

\def\TT{{\underline{T}}}
\def\vv{{\underline{v}}}
\def\tt{{\underline{t}}}

\def\dd{\underline{d}}

\def\mm{\underline{m}}

\def\1{\underline{1}}

\def\P{\Bbb P}

\def\Z{\Bbb Z}
\def\Q{\Bbb Q}
\def\C{\Bbb C}

\def\S{{\cal S}}
\def\O{{\cal O}}

\def\D{{\cal D}}
\def\Int{\mathsf{Int\,}}
\def\red{\mathsf{red\,}}

\newtheorem{theorem}{Theorem}

\newtheorem{lemma}{Lemma}

\newenvironment{corollary}
{\smallskip\noindent{\bf Corollary\/}.}{\smallskip\par}

\newenvironment{proof}
{\noindent{\bf Proof\/}.}{{ $\Box$}\smallskip\par}

\newenvironment{PROOF}
{\noindent{\bf Proof of Theorem 1\/}}{{ $\Box$}\smallskip\par}

\title{Universal abelian covers of rational surface singularities
and multi-index filtrations \footnote{Math. Subject Class. 14H20,
14J17, 32S25. Keywords: universal abelian covers, rational surface
singularities, Poincar\'e series.} }
\author{A.~Campillo
\and F.~Delgado \thanks{First two authors were partially supported
by the grant MTM2004-00958. Address: University of Valladolid, Dept.
of Algebra, Geometry and Topology, 47011 Valladolid, Spain. E-mail:
campillo\symbol{'100}agt.uva.es, fdelgado\symbol{'100}agt.uva.es}
\and S.M.~Gusein-Zade \thanks{Partially supported by the grants
RFBR-007-00593, INTAS-05-7805, NWO-RFBR 047.011.2004.026, and
RFBR-JSPS 06-01-91063. Address: Moscow State University, Faculty of
Mathematics and Mechanics, Moscow, GSP-2, 119992, Russia. E-mail:
sabir\symbol{'100}mccme.ru}}
\date{}
\begin{document}
\sloppy
\def\eps{\varepsilon}

\maketitle

\medskip

In \cite{inv} and \cite{coment}, there were computed the Poincar\'e
series of some (multi-index) filtrations on the ring of germs of
functions on a rational surface singularity. These Poincar\'e series
were written as the integer parts of certain fractional power
series, an interpretation of whom was not given. Here we show that,
up to a simple change of variables, these fractional power series
are specializations of the equivariant Poincar\'e series for
filtrations on the ring $\O_{\widetilde{\S},0}$ of germs of
functions on the universal abelian cover $(\widetilde{\S},0)$ of the
surface $(\S,0)$. We compute these equivariant Poincar\'e series.
From another point of view universal abelian covers of rational
surface singularities were studied in \cite{okuma}.

Let $(\S,0)$ be an isolated complex rational surface singularity
and let $\pi : (X,\D)\to (\S,0)$ be a resolution of it (not
necessarily the minimal one). Here $X$ is a smooth complex surface,
the exceptional divisor $\D = \pi^{-1}(0)$ is a normal crossing
divisor on $X$, all components $E_\sigma$ ($\sigma\in\Gamma$) of the
exceptional divisor $\D$ are isomorphic to the complex projective
line $\C\P^1$ and the dual graph of the resolution is a tree.

Let $\O_{\S,0}$ be the ring of germs of analytic functions on
$(\S,0)$. For $\sigma\in \Gamma$, i.e. for a component $E_\sigma$
of the exceptional divisor,
and for $f\in \O_{\S,0}$, let $v_\sigma(f)$ be the order of zero
of the lifting $f\circ \pi$ of the function $f$ to the space $X$
of the resolution along the component $E_\sigma$. Let us choose
several components $\seq E1s$ of the exceptional divisor $\D$
$(\{1,\ldots, s\}\subset \Gamma)$. The valuations $\seq v1s$
define a multi-index filtration $\{J(\vv)\}$
on the ring $\O_{\S,0}$:
for $\vv = (\seq v1s)\in \Z_{\ge 0}^s$,
$J(\vv) = \{f\in \O _{\S,0}\; : \; \vv(f)\ge \vv\}$ (here
$\vv(f) = (v_1(f), \ldots, v_s(f))\in \Z_{\geq 0}^s$,
$\vv'\geq \vv$ if and only if  $v'_i\geq v_i$ for all
$i=1,\ldots,s$). In \cite{inv}, there was computed the Poincar\'e
series $P (\seq t1s)$ of this filtration
(the definition of the Poincar\'e series of a multi-index
filtration can be found e.g. in \cite{inv, coment, mmj}). Let
$(E_\sigma\circ E_\delta)$ be the intersection matrix of the
components of the exceptional divisor. For $\sigma\neq \delta$,
the intersection number
$E_\sigma\circ E_\delta$ is equal to $1$ if the components
$E_\sigma$ and $E_\delta$ intersect (at one point) and is
equal
to zero if they don't intersect; the self-intersection number
$E_\sigma\circ E_\sigma$ of each component $E_\sigma$ is a
negative integer.
Let $d = \det (-(E_\sigma\circ E_\delta))$ and let
$(m_{\sigma \delta}) = - (E_\sigma\circ E_\delta)^{-1}$. All
entries $m_{\sigma \delta}$ are positive and
$m_{\sigma \delta}\in (1/d)\Z$. For $\sigma\in \Gamma$, let
$\mm_\sigma:= (m_{\sigma 1}, \ldots , m_{\sigma s})\in \Q_{\geq
0}^s$.

Let
$\stackrel{\bullet}{E_\sigma}$  be the ``smooth part" of the
component $E_\sigma$
in the exceptional divisor $\D$, i.e., $E_\sigma$ minus
intersection points with all other components of the
exceptional divisor ${\cal D}$.

For a fractional power series
$S(\seq t1s)\in \Z[[t_1^{1/d}, \ldots, t_{s}^{1/d}]]$, let
$\Int S(\seq t1s)$ be its ``integer part", i.e., the sum of all
monomials from $S(\seq t1s)$ with integer exponents.
In \cite{inv} it was shown that
\begin{equation}\label{eq1}
P(\seq t1s) = \Int \prod\limits_{\sigma\in\Gamma}
(1-\tt^{\mm_\sigma})^{-\chi(\stackrel{\bullet}{E_\sigma})}\; ,
\end{equation}
where $\tt^{\mm}:= t_1^{m_1}\cdot\ldots\cdot t_s^{m_s}$,
$\chi(X)$ is the Euler characteristic of the space $X$.

A similar formula was obtained in \cite{coment} for the Poincar\'e
series of the multi-index filtration on the ring $\O_{\S,0}$
defined by orders of a function germ on irreducible components of
a curve $(C,0) \subset (\S,0)$.

In \cite{inv}, the fractional power series
\begin{equation}\label{eq2}
Q(\tt) = \prod\limits_{\sigma\in
\Gamma}
(1-\tt^{\mm_\sigma})^{-\chi(\stackrel{\bullet}{E_\sigma})}
\end{equation}
(and
a similar one in \cite{coment}) participated as a formal
expression convenient to write the formula (\ref{eq1}) for the
Poincar\'e series $P(\seq t1s)$. There was no interpretation of it.

In \cite{mmj}, there was defined an equivariant Poincar\'e series
for an ``equivariant" filtration on the ring $\O_{V,0}$ of germs of
functions on a germ $(V,0)$ of a complex analytic variety with an
action of a finite group $G$. This Poincar\'e series was computed
for a divisorial filtration on the ring $\O_{\C^2,0}$ and for the
filtration defined by branches of a $G$-invariant plane curve
singularity $(C,0)\subset (\C^2,0)$ were the plane $\C^2$ was
equipped with a $G$-action.

Let $p: (\widetilde{\S},0)\to (\S,0)$ be the universal abelian cover
of the surface singularity $(\S,0)$: see e.g. \cite{NW, okuma}. One
can describe it in the following way. Let G $=
H_1(\S\setminus\{0\})$ be the first homology group of the
(nonsingular) surface $\S\setminus\{0\}$. The order of the group $G$
is equal to the determinant $d$ of the minus intersection matrix
$-(E_\sigma\circ E_\delta)$ and moreover  $G$ is the cokernel
$\Z^\Gamma/ \mbox{Im\,}I$ of the map $I: \Z^{\Gamma}\to \Z^{\Gamma}$
defined by this matrix.

The group $G$ acts on the germ $(\widetilde{\S},0)$ and the
restriction $p|_{\widetilde{\S}\setminus\{0\}}$ of the map $p$ to
the complement of the origin is a (usual, nonramified) covering
$\widetilde{\S}\setminus \{0\} \to \S\setminus \{0\}$  with the
structure group $G$. One can lift the map $p$ to a (ramified)
covering $p: (\widetilde{X}, \widetilde{\D})\to(X,\D)$ where
$\widetilde{X}$ is a normal surface (generally speaking not smooth)
and $\widetilde{X}\setminus \widetilde{\D} \cong
\widetilde{\S}\setminus \{0\}$:
$$
\begin{array}{ccc}
(\widetilde{X}, \widetilde{\D}) &
\stackrel{\widetilde{\pi}}{\longrightarrow}
& (\widetilde{\S},0) \\
\  \downarrow {\small p} &  & {\small p} \downarrow \ \\
(X, \D) & \stackrel{\pi}{\longrightarrow}  & (\S,0) \\
\end{array}
$$
(one can define $\w{X}$ as the normalization of the fibre
product $X\times_{\S}\widetilde{\S}$ of the varieties $X$ and
$\widetilde{\S}$ over $\S$).

Let $g_\sigma$, $\sigma \in \Gamma$ be the element of the group $G$
represented by the loop in $X\setminus{\D}$ going around the
component $E_\sigma$ in the positive direction. The group $G$ is
generated by the elements $g_\sigma$ for all $\sigma\in\Gamma$. For
a point $x\in \stackrel{\bullet}{E_\sigma}$ and for a point
$\widetilde{x}$ from the preimage $p^{-1}(x)$ of it, locally, in a
neighbourhood of the point $\widetilde{x}$, the map $p:
\widetilde{D}\to \D$ is an isomorphism and the map $p:
\widetilde{X}\to X$ is a ramified (over $\D$) covering, the order
$d_\sigma$ of which coincides with the order of the generator
$g_\sigma$ of the group $G$.

\begin{lemma}\label{lemma1}
The order $d_\sigma$ of the element $g_\sigma\in G$ is
the minimal natural $k$ such that $k m_{\delta \sigma}$ is an
integer for all $\delta\in \Gamma$.
\end{lemma}

\begin{proof}
This follows immediately from the fact that $\Z^{\Gamma}/\mbox{Im\,}
I\cong \mbox{Im\,} m/ \Z^{\Gamma}$ where $m: \Z^{\Gamma}\to
\Q^{\Gamma}$ is the map given by the matrix $(m_{\sigma \delta})$
(i.e. minus the inverse of the map $I$).
\end{proof}

Let $R(G)$ be the ring of (virtual) representations of the
group $G$.
For $\sigma\in \Gamma$, let $\alpha_\sigma$ be the
one-dimensional representation  $G\to \C^* = {\bf GL}(1,\C)$
of the group $G$ defined by
$\alpha_\sigma(g_\delta) = \exp (- 2\pi \sqrt{-1} m_{\sigma
\delta})$
(here the minus sign reflects the fact that the action of an
element $g\in G$ on the ring $\O_{\widetilde{\S},0}$ is defined by
$(g\cdot f)(x) = f(g^{-1}(x))$).

Let us choose any component $\widetilde{E}_i$ of the preimage
$p^{-1}(E_i)$ of the component $E_i$ and let $\widetilde{v}_i$ be
the corresponding divisorial valuation on the ring
$\O_{\widetilde{\S},0}$. On the space $\bigcup\limits_{\alpha}
\O_{\S,0}^{\alpha}$ of all $G$-equivariant functions on $(\w{\S},0)$
($\alpha$ runs over all nonequivalent 1-dimensional representations
of the group $G$) the valuation $\w{v}_i$ does not depend on the
choice of the component $\w{E}_i$.

In \cite{mmj}, there was defined the equivariant Poincar\'e series
of the multi-index filtration defined by the divisorial valuations
$\seq{\w{v}}1s$.

\begin{theorem}\label{theo1}
The equivariant Poincar\'e series $P^{G}(\seq t1s)$ of the
$s$-index filtration defined by the set of
divisorial valuations $\{\seq{\widetilde{v}}{1}{s}\}$ is given by
the formula:
\begin{eqnarray}
\label{eq3}
P^G(\seq t1s) &=&
\prod_{\sigma\in \Gamma} (1- \alpha_\sigma
\tt^{\dd\, \mm_\sigma})^{-\chi(\stackrel{\bullet}{E_\sigma})}
\nonumber
\\
&=&
\prod_{\sigma\in \Gamma} (1- \alpha_\sigma
t_1^{d_1 m_{1 \sigma}}
\cdot\ldots\cdot
t_s^{d_s m_{s \sigma}})^{-\chi(\stackrel{\bullet}{E_\sigma})} \;
.
\end{eqnarray}
\end{theorem}

For a power series $S(\seq t1s) = \sum\limits_{\vv\in
\Z^s_{\geq 0}} s_{\vv}\tt^{\vv}\in R(G)[[\seq t1s]]$ ($R(G)$ is
the ring of representations of the group $G$),  let its
reduction
$\red S(\seq t1s)$ be the series
$\sum\limits_{\vv\in
\Z^s_{\geq 0}} (\dim s_{\vv})\tt^{\vv}\in \Z [[\seq t1s]]$.

\begin{corollary}
One has $\red P^G(t_1, \ldots, t_s) = Q(t_1^{d_1}, \ldots,
t_s^{d_s})$, where $Q(\tt)$ is the fractional power series
defined by (\ref{eq2}).
\end{corollary}

\begin{PROOF}
For short we shall say that an effective divisor on
$\stackrel{\bullet}{\D} = \bigcup \stackrel{\bullet}{E_\sigma}$ (or
on $\stackrel{\bullet}{\widetilde{\D}} =
p^{-1}(\stackrel{\bullet}{\D})$) is Cartier if it is the
intersection with  $\stackrel{\bullet}{\D}$ (or with
$\stackrel{\bullet}{\widetilde{\D}}$) of the strict transform of a
Cartier divisor on $(\S,0)$ (or on $(\widetilde{\S},0)$). From
\cite{mmj} it follows that the equivariant Poincar\'e series
$P^G(\tt)$ is equal to the integral with respect to the Euler
characteristic of the monomial $\alpha\,\tt^{\widetilde{\vv}}$ over
the space of $G$-invariant effective Cartier divisors on
$\stackrel{\bullet}{\widetilde{\D}}$. Here $\alpha\in R(G)$ and
$\widetilde{\vv}\in \Z^s_{\geq 0}$ are functions (in fact semigroup
homomorphisms) on the space of $G$-invariant Cartier divisors on
$\stackrel{\bullet}{\widetilde{\D}}$: a $G$-invariant Cartier
divisor defines the orders of zero of the corresponding
($G$-equivariant) function along the components
$\stackrel{\bullet}{E_i}$ and also the corresponding 1-dimensional
representation of the group $G$.

Thus to compute the equivariant Poincar\'e series $P^G(\tt)$ one
has to describe the space of $G$-invariant effective Cartier
divisors on  $\stackrel{\bullet}{\widetilde{\D}}$ and the
corresponding functions $\vv$ and $\alpha$ on it.

\begin{lemma}\label{lemma2}
Any $G$-invariant effective divisor on
$\stackrel{\bullet}{\widetilde{\D}}$ is a Cartier divisor.
\end{lemma}

\begin{proof}
It is sufficient to show this for the divisor
$\sum\limits_{\widetilde{x}\in p^{-1}(x)} \widetilde{x}$ for a point
$x\in \stackrel{\bullet}{E_\sigma}$, i.e. for the $G$-orbit of a
point from $\stackrel{\bullet}{\widetilde{E}_\sigma}$. The isotropy
group $G_{\widetilde{x}}$ of a point $\widetilde{x}\in p^{-1}(x)$ is
the cyclic subgroup of the group $G$ of order $d_\sigma$ generated
by the element $g_\sigma$ (this element acts trivially on
$p^{-1}(E_\sigma)$).

Let us take the germ at the point $x$ of
a smooth curve $L_\sigma$ on $(X,\D)$
transversal to the exceptional divisor $\D$.
By the
Artin criterion (see, e.g., \cite{P}, Lemma on page 156), the
divisor $d\cdot L_\sigma$ is the strict transform of a Cartier
divisor
on $(\S,0)$ (in fact already $d_\sigma L_\sigma$ is one with
this property), i.e. there exists a function $f_\sigma: \S\to\C$
such that the strict
transform of the divisor $\{f_\sigma=0\}$ is $d\cdot L_\sigma$.
Let
$\overline{f}_\sigma = f_\sigma \circ \pi$ be the lifting of the
function $f_\sigma$ to the space $X$ of the resolution and let
$\widetilde{f}_\sigma = f_\sigma\circ \pi \circ p$ be
the lifting of the function $f_\sigma$ to the space
$\widetilde{X}$ of the
modification of the universal abelian cover $(\widetilde{\S},0)$
($\widetilde{f}_\sigma$ is a $G$-invariant function on
$\widetilde{X}$).
Let us describe the divisor $\{\widetilde{f}_\sigma =0\}$.
Let $\widetilde{L}_{\sigma, \widetilde{x}}\subset \widetilde{X}$
be  the germ at the point $\widetilde{x}\in p^{-1}(x)$ of the
preimage under the map $p$ of the curve $L_\sigma\subset X$.

The order of zero of the function $\widetilde{f}_\sigma$ along
$\widetilde{L}_{\sigma, \widetilde{x}}$ is equal to $d$.  The order
of zero of the function $\overline{f}_\sigma$ along the component
$\stackrel{\bullet}{E}_{\delta}$ is equal to $d\cdot m_{\sigma
\delta}$. The ramification order of the map $p$ over the component
$\stackrel{\bullet}{E}_{\delta}$ is equal to $d_\delta$. Therefore
the order of zero of the function $\widetilde{f}_\sigma =
\overline{f}_\sigma\circ p$ along the preimage of the component
$\stackrel{\bullet}{E}_{\delta}$ is equal to $d\cdot d_\delta\cdot
m_{\sigma \delta}$. This (integer) number is divisible by $d$ (since
$d_\delta m_{\sigma \delta}$ is an integer: see Lemma~\ref{lemma1}).
Therefore the zero divisor of the function $\widetilde{f}_\sigma$ is
divisible by $d$, i.e. the order of zero of this function along each
component of its zero set is divisible by $d$. This means that a
root $\sqrt[d]{\widetilde{f}_\sigma}$ of degree $d$ of the function
$\widetilde{f}_\sigma$ (i.e. a branch of this root) is well defined
up to multiplication by a root of degree $d$ of $1$ $G$-equivariant
complex analytic function on $\widetilde{X}$ and thus it is the
lifting  of a $G$-equivariant function on $(\widetilde{\S},0)$ (see
e.g. \cite[page ?]{lau}).
\end{proof}

\begin{corollary}
Each $G$-invariant divisor on the universal abelian cover
$(\widetilde{\S},0)$ of the rational surface singularity $(\S,0)$
is a Cartier one.
\end{corollary}

Lemma \ref{lemma2} means that the space of $G$-invariant
effective Cartier divisors on
$\stackrel{\bullet}{\widetilde{\D}}$  is in one to one
correspondence with the space of all effective divisors on
$\stackrel{\bullet}{\widetilde{\D}}$. As it follows from the
proof of Lemma \ref{lemma2}, the order of zero of the
$G$-equivariant
function $\widetilde{f}_\sigma$ (corresponding to one point
${x}\in \stackrel{\bullet}{E_\sigma}$) along the
component $\widetilde{E}_i$ is equal to $d_i m_{\sigma i}$.
One has to find the (one-dimensional) representation
$\alpha_\sigma$ with respect to which the function
$\widetilde{f}_\sigma$ is $G$-equivariant.

\begin{lemma}
$$
\alpha_\sigma (g_{\delta}) = \exp (- 2 \pi \sqrt{-1} m_{\sigma
\delta})\; .
$$
\end{lemma}

\begin{proof}
The element $g_\delta$ of the group $G$ acts trivially on the
preimage
$p^{-1}(\stackrel{\bullet}{E_\delta})$ of the component
$\stackrel{\bullet}{E_\delta}$ of the exceptional divisor and
acts by multiplication by
$\exp(\dfrac{2 \pi}{d_\delta} \sqrt{-1})$ on the normal line to
it. The order of zero of the function $\widetilde{f}_\sigma$
along the preimage
$p^{-1}(\stackrel{\bullet}{E_\delta})$ is equal to
$m_{\sigma  \delta} d_\delta$. Therefore
$$
\frac{g_{\delta}\cdot f_\sigma}{f_\sigma} = \exp(- \frac{ 2
\pi\sqrt{-1} m_{\sigma \delta} d_\delta }{d_\delta})
= \exp(-  2
\pi\sqrt{-1} m_{\sigma \delta})\; .
$$
\end{proof}

Now Theorem follows from the usual arguments used e.g. in \cite{inv,
mmj}. The space of effective divisors on $\stackrel{\bullet}{\D} =
\bigcup \stackrel{\bullet}{E}_\sigma$ is the direct product of the
spaces of effective divisors on the components
$\stackrel{\bullet}{E}_\sigma$. Each of the latter ones is the
disjoint union of symmetric  powers $S^k
\stackrel{\bullet}{E}_\sigma$ of the component
$\stackrel{\bullet}{E}_\sigma$. Therefore
$$
P^G (\seq t1s) = \prod_{\sigma\in \Gamma}
\left(\sum_{k=0}^{\infty}
\chi (S^k \stackrel{\bullet}{E}_\sigma) \cdot \alpha_\sigma^k \tt^{k
\dd\,\mm_\sigma} \right)\; ,
$$
(this follows from the fact that $\vv$ and $\alpha$ are
semigroup homomorphisms).
The well-known formula
$$
\sum_{k=0}^{\infty}
\chi (S^k X) t^k = (1-t)^{-\chi(X)}
$$
implies the equation~(\ref{eq3}).
\end{PROOF}

A similar result holds for the filtration on the ring
$\O_{\w{\S},0}$
defined
by orders of a function germ on branches of a $G$-invariant curve
$(\widetilde{C},0)\subset (\widetilde{\S},0)$.
Let $\widetilde{C} = \bigcup\limits_{i=1}^r \widetilde{C}_i$
where
$\widetilde{C}_i$ are irreducible $G$-invariant components of the
curve
$\widetilde{C}$ (generally speaking each curve $\widetilde{C}_i$
consists of several irreducible compoents permuted by the group
$G$). Each curve
$\widetilde{C}_i$ is the preimage under the map $p$ of an
irreducible curve  $C_i$ on $(\S,0)$.
The curve $\w{C}=\bigcup\limits_{i=1}^r\w{C}_i$ defines an
$r$-index filtration on the space
$\bigcup\limits_{\alpha} \O_{\w{\S},0}^\alpha$ of $G$-equivariant
functions on the surface $(\w{\S},0)$ (or on  the space
$\bigcup\limits_{\alpha} \O_{\w{C},0}^\alpha$ of $G$-equivariant
functions on the cuvre $(\w{C},0)$).
Let $\varphi_i : (\C,0)\to (\w{\S},0)$ be a parametrization
(uniformization) of an irreducible component of the curve
$\w{C}_i$. For a $G$-equivariant function germ $f$, let
$\w{w}_i(f)$ be the order of zero of the function
$f\circ \varphi_i$ at the origin: $f\circ \varphi_i(\tau)= a
\tau^{\w{w}_i(f)} +$ terms of higher degree, $a\neq 0$.
The valuations $\seq{\w{w}}1r$  define a multi-index filtration
in the usual way.

Let $\pi: (X,\D)\to (\S,0)$ be a resolution of the surface
singularity $(\S,0)$ which at the same time is an embedded
resolution of the curve $(C,0)\subset (\S,0)$, $C =
\bigcup\limits_{i=1}^r C_i$. Let $\overline{C}_i$ be the strict
trasnform of the curve $C_i$ in $X$. Let $\seq E1s$ be all the
components of the exceptional divisor $\D$ of the resolution. Let
$\stackrel{\circ}{E}_i$ be the ``smooth part" of the component $E_i$
in the total transform $\pi^{-1}(C)$ of the curve $C$, i.e. $E_i$
minus intersection points with all other components of the total
transform $\pi^{-1}(C)$. Let $\mm_i= (m_{i1}, \ldots, m_{is})\in
\Q_{\geq 0}^{s}$, $\dd = (\seq{d}1s)\in \Z^s_{\geq 0}$, and a
1-dimensional representation $\alpha$ of the group $G$ ($i=1,\ldots,
s$) be defined as above. The same arguments as in the proof of
Theorem 1 imply the following statement.

\begin{theorem}\label{theo2}
The equivariant Poincar\'e series $P^{G}(\seq t1r)$ of the
$r$-index filtration defined by the set of
valuations $\{\seq{\widetilde{w}}{1}{r}\}$ is given by
the formula:
$$
P^G(\seq t1r) = \left( \prod_{i=1}^s (1- \alpha_i \TT^{\dd\,
\mm_i})^{-\chi(\stackrel{\circ}{E_i})} \right)
\mbox{\raisebox{-2ex}{ $\vert_{\; T_i\; \mapsto \prod\limits_{j \,
:\, \overline{C}_j\cap E_i = pt} t_j}$}} \;
$$
$($here $\TT = (\seq{T}1s)$; in the substitution above,
$\prod\limits_{j\in \emptyset}t_j$ is supposed to be equal to
$1$$)$.
\end{theorem}

\end{document}